\documentclass{amsart}

\usepackage{amsmath,amssymb,amsthm,latexsym}

\newtheorem{theorem}{Theorem}
\newcommand{\bt}{\begin{theorem}}
\newcommand{\et}{\end{theorem}}
\newtheorem{lemma}{Lemma}
\newcommand{\bl}{\begin{lemma}}
\newcommand{\el}{\end{lemma}}
\newtheorem{corollary}{Corollary}
\newcommand{\bc}{\begin{corollary}}
\newcommand{\ec}{\end{corollary}}
\newcommand{\bconj}{\begin{conjecture}}
\newcommand{\econj}{\end{conjecture}}
\newtheorem{problem}{Problem}
\newcommand{\bprob}{\begin{problem}}
\newcommand{\eprob}{\end{problem}}
\newcommand{\beq}{\begin{equation}}
\newcommand{\eeq}{\end{equation}}
\newcommand{\benum}{\begin{enumerate}}
\newcommand{\eenum}{\end{enumerate}}
\newcommand{\N}{\ensuremath{ \mathbf N }}
\newcommand{\Z}{\ensuremath{\mathbf Z}}
\newcommand{\Q}{\ensuremath{\mathbf Q}}
\newcommand{\R}{\ensuremath{\mathbf R}}

\newcommand{\mcx}{\ensuremath{ \mathcal X}}

\newcommand{\mba}{\ensuremath{ \mathbf a}}

\newcommand{\mbo}{\ensuremath{ \mathbf 0}}

\newcommand{\mbx}{\ensuremath{ \mathbf x}}
\newcommand{\mby}{\ensuremath{ \mathbf y}}

\newcommand{\bmat}{\left(\begin{matrix}}
\newcommand{\emat}{\end{matrix}\right)}
\newcommand{\bsmallmat}{\left(\begin{smallmatrix}}
\newcommand{\esmallmat}{\end{smallmatrix}\right)}

\DeclareMathOperator{\qqand}{\qquad\text{and}\qquad}

\title[${\mathbf Q}$-independence and $B_h$-sets]{${\mathbf Q}$-independence and the construction  of $B_h$-sets of integers and lattice points} 
\author{Melvyn B.  Nathanson}
\address{Department of Mathematics\\Lehman College (CUNY)\\Bronx, NY 10468}
\email{melvyn.nathanson@lehman.cuny.edu}

\date{\today}

\thanks{Supported in part by  PSC-CUNY Research Award Program grant 66197-00 54.}

\begin{document}

\begin{abstract}
This paper gives a simple ${\mathbf Q}$-vector space construction  of  finite $B_h$-sets of integers and lattice points.
\end{abstract}

\maketitle

\section{$B_h$-subsets of a group} 
Let $G$ be an additive abelian group or semigroup and let $A$ be a nonempty subset of $G$.  
For every positive integer $h$, the $h$-fold \emph{sumset}  of $A$  is the set of all sums of $h$ not necessarily distinct elements of $A$:
\[
hA = \underbrace{A+\cdots + A}_{\text{$h$ summands}}
= \left\{ a'_1+\cdots + a'_h: a'_i \in A \text{ for all } i \in [1,h] \right\}. 
\]
For every element $g \in G$, 
the representation function $r_{A,h}(g)$ counts the number of 
representations of $g$ as a sum of $h$ not necessarily distinct elements 
of $A$.  Representations that differ only by a permutation of the summands 
are not counted separately.  
The set $A$ is a \emph{$B_h$-set} if $r_{A,h}(g) = 0$ or $1$ for all  $g \in G$.  
Thus, a $B_h$-set is a set $A$ such that the representation of an element of the 
semigroup as the sum of $h$ elements of $A$ is unique, up to permutation of the summands.  
$B_2$-sets are also called \emph{Sidon sets}.

The study of $B_h$-sets is a classical topic in combinatorial additive number theory
(O'Bryant~\cite{obry04}).  Although ``almost all'' finite sets of integers are $B_h$-sets 
(Nathanson~\cite{nath2004-105}) and there are many papers estimating the size such sets, 
there are few explicit constructions of them 
(cf.~\cite{bose-chow62}--\cite{ruzs93a}).
This paper gives a simple \Q-vector space construction of finite $B_h$-sets 
for the group \Z\ of integers and the group $\Z^d$ of integer lattice points.

Let $\N = \{1,2,3,\ldots\}$ be the set of positive integers and 
 $\N_0 = \N \cup \{0\} = \{0,1,2,3\ldots\}$  the set of nonnegative integers.  
For $u,v \in \R$, the \emph{interval of integers} $[u,v]$ 
is  the set of all integers $k$ such that $u \leq k \leq v$. 
The \emph{integer part} of the real number $u$ (also called the \emph{floor} of $u$), 
is the largest integer $n \leq u$.

\section{The set $\mcx_{h,n}$}
For positive integers $h$ and $n$, let  
\[
\mcx_{h,n} = \left\{ \mbx = (x_1,\ldots, x_n) \in \N_0^n: \sum_{i=1}^n x_i = h \right\}. 
\] 
We have  $\left|  \mcx_{h,n} \right| = \binom{n+h-1}{h}$.
The set $\mcx_{h,n}$ is symmetric in the sense that if $\mbx = (x_1,\ldots, x_h) \in \mcx_{h,n}$, 
then $\sigma \mbx = \left(x_{\sigma(1)}, \ldots,x_{\sigma(n)} \right) \in \mcx_{h,n}$ 
for every permutation $\sigma$ of $[1,n]$.

Let $A = \{a_1,\ldots, a_n\}$ be a  finite subset of an additive abelian 
semigroup with $|A| = n$.   
The $h$-fold sumset of  $A $ can  be written in the form  
\beq                \label{Sidon:hA-Xform}
hA = \left\{ \sum_{i=1}^n x_i a_i : \mbx   = (x_1,\ldots, x_n) \in \mcx_{h,n} \right\}. 
\eeq

For all $\mbx = (x_1,\ldots, x_n) \in \mcx_{h,n}$ and $\mby = (y_1,\ldots, y_n) \in \mcx_{h,n}$
with $\mbx \neq \mby$, we have  
\beq                \label{Sidon:xy-sum}
2 \leq \sum_{i=1}^n |x_i - y_i| \leq \sum_{i=1}^n  x_i  +  \sum_{i=1}^n  y_i = 2h. 
\eeq
These upper and lower bounds are best possible.  
If $n \geq h$ and 
\[
\mbx =  (2,0,\underbrace{1,1,\ldots, 1}_{\text{$(h-2)$ 1's}}, \underbrace{0,0,\ldots, 0}_{\text{$(n-h)$ 0's}} )
\] 
and 
\[
\mby =  (1,1,\underbrace{1,\ldots, 1}_{\text{$(h-2)$ 1's}}, \underbrace{0,0,\ldots, 0}_{\text{$(n-h)$ 0's}} )
\]
then $ \sum_{i=1}^n |x_i - y_i| = 2$.

If $n \geq 2h$ and 
\[
\mbx =  (\underbrace{1,1,\ldots, 1}_{\text{$h$ 1's}}, 
\underbrace{0,0,\ldots, 0}_{\text{$h$ 0's}}\underbrace{0,0,\ldots, 0}_{\text{$(n-2h)$ 0's}} )
\] 
and 
\[
\mby =  (\underbrace{0,0,,\ldots, 0}_{\text{$h$ 1's}},\underbrace{1,1,\ldots, 1}_{\text{$h$ 1's}}, \underbrace{0,0,\ldots, 0}_{\text{$(n-2h)$ 0's}}  )
\]
then $ \sum_{i=1}^n |x_i - y_i| = 2h$.

\section{$B_h$-sets of lattice points}

The $\ell^{\infty}$-norm of the vector $\vec{\mba} = (a_1, a_2, \ldots, a_d) \in \R^d$ is
\[
\|\vec{\mba} \|_{\infty} = \max_{j \in [1,d]} |a_j|.
\]
The $\ell^{\infty}$-norm of the set of vectors $A = \{ \vec{\mba}_1, \ldots, \vec{\mba}_n \}$ is
\[
\| A \|_{\infty} = \max_{i \in [1,n]} \|\vec{\mba}_i \|_{\infty}.
\]

Let $  \left( \vec{\theta}_1, \vec{\theta}_2, \ldots, \vec{\theta}_n \right)$
 be an $n$-tuple of \Q-independent vectors in $\R^d$ 
 and let 
 \[
\Theta_n = \left\{\vec{\theta}_1, \vec{\theta}_2, \ldots, \vec{\theta}_n \right\} 
 \] 
 and 
 \[
\| \Theta_n \|_{\infty} = \max_{i \in [1,n]} \|\vec{\theta}_i \|_{\infty}.
\]
For all $i \in [1,n]$, we write 
\[
\vec{\theta}_i = (\theta_{i,1}, \theta_{i,2},\ldots, \theta_{i,d}).
\]
The \Q-independence of the set $\Theta_n$ implies that, 
for all $\vec{\mbx}, \vec{\mby} \in \mcx_{h,n}$ 
with $\mbx = (x_1,\ldots, x_n)$, $\mby = (y_1,\ldots, y_n)$, 
and $\vec{\mbx} \neq \vec{\mby}$, we have 
\[
\sum_{i=1}^h x_i \vec{\theta}_i \neq \sum_{i=1}^h y_i \vec{\theta}_i.
\] 
Equivalently, 
\[
  \sum_{i=1}^n (x_i -  y_i) \vec{\theta}_i \neq \mbo.  
\]
Inequality~\eqref{Sidon:xy-sum} gives   
\begin{align*}
0 & <  \left\|   \sum_{i=1}^n (x_i -  y_i) \vec{\theta}_i \right\|_{\infty}  
 \leq \| \Theta_n \|_{\infty}  \sum_{i=1}^n |x_i -  y_i|  \leq  2h\| \Theta_n \|_{\infty}. 
\end{align*} 
Because the set $\mcx_{h,n}$ is finite, we have 
\beq                  \label{Sidon:varepsilon}
0 < \varepsilon_{h,n}  = \inf \left\{ \left\|   \sum_{i=1}^n (x_i -  y_i) \vec{\theta}_i \right\|_{\infty} :
\mbx,\mby \in \mcx_{h,n} \text{ and } \mbx \neq \mby  \right\}.  
%\leq 2h\| \Theta_n \|_{\infty}.  
\eeq

Let $q$ and $m$ be positive integers.  For all $i \in [1,n]$ and $j \in [1,d]$, 
there are exactly $2m$ integers $a_{i,j}$ 
such that 
\beq                       \label{Sidon:aq-theta}
0 < |a_{i,j}-q\theta_{i,j}| \leq m.
\eeq 
For every integer  $a_{i,j}$  such that $|a_{i,j}-q\theta_{i,j}| \leq m$, we have  
\[     
\left| \frac{a_{i,j}}{q} -  \theta_{i,j} \right| \leq \frac{m}{q}. 
\] 
Note that  if  $\theta_{i,j} \geq 0$, then there are at least $m$ positive integers $a_{i,j}$ 
that satisfy inequality~\eqref{Sidon:aq-theta}.     

For all $i \in [1,n]$, we obtain $(2m)^d$ lattice points    
\[
\vec{a}_i = (a_{i,1}, a_{i,2}, \ldots, a_{i,d}) \in \Z^d 
\]
and $(2m)^d$ rational vectors  
\[
\frac{1}{q} \vec{a}_i = \left(  \frac{a_{i,1}}{q},  \frac{a_{i,2}}{q}, \ldots,, \frac{a_{i,d}}{q}\right) \in \Q^d.
\]
These vectors approximate  $\vec{\theta}_i$:  
\begin{align}
\left\| \frac{1}{q} \vec{a}_i  - \vec{\theta}_i\right\|_{\infty}  
& = \left\| \left( \frac{a_{i,1}}{q} -  \theta_{i,1},  \frac{a_{i,2}}{q} -  \theta_{i,2}, \ldots, 
 \frac{a_{i,d}}{q} -  \theta_{i,d}\right) \right\|_{\infty}        \notag \\
 & = \max_{j \in [1,d]} \left| \frac{a_{i,j}}{q} -  \theta_{i,j} \right|     \label{Sidon:1/q-ineq}  \\
 &  \leq \frac{m}{q}.     \notag  
\end{align} 
Thus, we have constructed $(2m)^{dn}$ sets of lattice points  
\[
A_{h,n}(q,m) = \{  \vec{a}_1,\ldots,  \vec{a}_n\} \subseteq \Z^d.  
\]
These sets are  determined by the positive integers $q$ and $m$ and by the choices of the  
$dn$ coordinates $a_{i,j}$.

If $\vec{\theta}_i$ is a nonnegative vector, then we can choose  coordinates  $a_{i,j}$ 
so that the  lattice point $\vec{\mba}_i$ has positive coordinates.  
If $\Theta$ is a set of nonnegative vectors, 
then we can choose  positive  coordinates $a_{i,j}$ for all $n$  vectors  
in the set $A_{h,n}(q,m)$.

\bt              \label{Sidon:theorem:construction}
Let $\Theta_n = \{\vec{\theta}_1, \vec{\theta}_2, \ldots, \vec{\theta}_n \}$ 
be a set of \Q-independent vectors in $\R^d$.  
Let $h \geq 2$ and define $\varepsilon_{h,n}$ by~\eqref{Sidon:varepsilon}.   
For all positive integers $q$ and $m$ with  
\[
q > \frac{2hm}{\varepsilon_{h,n}} 
\]
the  $(2m)^{dn}$ sets $A_{h,n}(q,m)$ 
 constructed from $\Theta_n$  are $B_h$-sets of lattice points in $\Z^d$ with 
\beq                  \label{Sidon:A-Theta-norm} 
\| A_{h,m}(q) \|_{\infty} \leq q \| \Theta_n \|_{\infty}+ m. 
\eeq
\et

\begin{proof} 
Let  $A_{h,n}(q,m) = \{\vec{a}_1, \vec{a}_2,\ldots, \vec{a}_n\}$.  
Inequality~\eqref{Sidon:aq-theta} implies that, for all $i \in [1,n]$ and $j \in [1,d]$, 
\[
|a_{i,j}| \leq q|\theta_{i,j}| + m  \leq q \left\|\vec{\theta}_i \right\|_{\infty} + m
\]
and so 
\[
\| \vec{a}_i \|_{\infty} \leq q \|\vec{\theta}_i\|_{\infty} + m \leq q \| \Theta_n \|_{\infty}+ m. 
\]
This immediately gives the upper bound~\eqref{Sidon:A-Theta-norm}.

We shall prove that $A_{h,n}(q,m)$ is a $B_h$-set of lattice points in $\Z^n$.  
By formula~\eqref{Sidon:hA-Xform}, we have the $h$-fold sumset 
\[
hA_{h,n}(q,m) = \left\{  \sum_{i=1}^n x_i \vec{a}_i : \mbx = (x_1,\ldots, x_d) \in \mcx_{h,n} \right\}.
\]
Let $\vec{\mbx} = (x_1,\ldots, x_n) \in \mcx_{h,n}$ and $\vec{\mby} = (y_1,\ldots, y_n) \in \mcx_{h,n}$ 
with $\mbx \neq \mby$.  Applying inequalities~\eqref{Sidon:xy-sum} 
and~\eqref{Sidon:1/q-ineq},  we obtain  
 \begin{align*}
\left\|  \sum_{i=1}^n (x_i-y_i)\left( \frac{1}{q} \vec{a}_i  - \vec{\theta}_i \right) \right\|_{\infty}
& \leq \sum_{i=1}^n |x_i-y_i|    \left\|  \frac{1}{q} \vec{a}_i  - \vec{\theta}_i \right\|_{\infty} \\
 &  \leq \frac{m}{q}\sum_{i=1}^n |x_i-y_i|   \\
 & \leq  \frac{2hm}{q}  
\end{align*}
and so   
 \begin{align*}
\frac{1}{q} & \left\|   \sum_{i=1}^n x_i \vec{a}_i -  \sum_{i=1}^n y_i \vec{a}_i \right\|_{\infty}  
 = \left\|   \sum_{i=1}^n (x_i - y_i) \frac{1}{q} \vec{a}_i \right\|_{\infty}  \\
& = \left\|   \sum_{i=1}^n (x_i - y_i)  \vec{\theta}_i 
+ \sum_{i=1}^n (x_i - y_i) \left( \frac{1}{q} \vec{a}_i  - \vec{\theta}_i \right)  \right\|_{\infty}  \\
& \geq \left\|   \sum_{i=1}^n (x_i - y_i)  \vec{\theta}_i  \right\|_{\infty}   - 
\left\|  \sum_{i=1}^n (x_i - y_i) \left( \frac{1}{q} \vec{a}_i  - \vec{\theta}_i \right)  \right\|_{\infty}  \\ 
& \geq \varepsilon_{h,n} -  \frac{2hm}{q}. 
\end{align*} 
Choosing an integer 
\[
q > \frac{2hm}{\varepsilon_{h,n}}
\]
we obtain 
\[
 \left\|   \sum_{i=1}^n x_i \vec{a}_i -  \sum_{i=1}^n y_i \vec{a}_i \right\|_{\infty}  
 \geq q\varepsilon_{h,n} - 2hm > 0
\]
and so 
$A_{h,n}(q,m) = \{\vec{a}_1,\ldots, \vec{a}_n\}$ is a $B_h$-set of lattice points.  
This completes the proof.  
\end{proof}

Applying Theorem~\ref{Sidon:theorem:construction} in dimension $d=1$ 
gives the following result for $B_h$-sets of integers.

\bt              \label{Sidon:theorem:construction-integers}
Let $\Theta_n = \{ \theta_1, \theta_2, \ldots, \theta_n \}$ 
be a set of \Q-independent real numbers.
Let $h \geq 2$ and define $\varepsilon_{h,n}$ by~\eqref{Sidon:varepsilon}.   
For all positive integers $q$ and $m$ with  
\[
q > \frac{2hm}{\varepsilon_{h,n}} 
\]
the  $(2m)^{n}$ sets $A_{h,n}(q,m)$ 
 constructed from $\Theta_n$  are$B_h$-sets of integers with 
\[
\| A_{h,n}(q) \|_{\infty} \leq q \| \Theta_n \|_{\infty}+ m. 
\]
\et

\section{Examples of the construction of Sidon sets}

A Sidon set is a $B_2$-set.  
The set $A$ of size $n$ is a Sidon set if and only if $|2A| = \binom{n+1}{2}$. 
For all  integers $n \geq 2$, 
\[
\mcx_{2,n} = \left\{ \mbx = (x_1,\ldots, x_n) \in \N_0^n: \sum_{i=1}^n x_i = 2 \right\} 
\] 
and $|\mcx_{2,n} | = \binom{n+1}{2}$. 
We have 
\begin{align*}
\mcx_{2,1} & =  \left\{ (2)  \right\} \\
\mcx_{2,2} & =  \left\{ (2,0), (1,1), (0,2) \right\} \\
\mcx_{2,3} & =  \left\{ (2,0,0), (0,2,0), (0,0,2),  (1,1,0), (1,0,1), (0,1,1) \right\} \\
\mcx_{2,4} & =   \{ (2,0,0,0), (0,2,0,0), (0,0,2,0), (0,0,0,2),    (1,1,0,0),  \\
& \qquad (1,0,1,0), (1,0,0,1 ),  (0,1,1,0), (0,1,0,1), (0,0,1,1) \} 
\end{align*}
Let $\Theta_n = \left\{\theta_1, \ldots, \theta_n \right\}$ 
be a \Q-independent set of $n$ real numbers.
Then  
\[
\varepsilon_{2,n}(\Theta_n) 
= \inf \left\{ \left|   \sum_{i=1}^n (x_i -  y_i) \theta_i \right| :
\mbx,\mby \in \mcx_{2,n} \text{ and } \mbx \neq \mby  \right\} 
\]
and so 
\begin{align*}
\varepsilon_{2,2}(\Theta_2) 
& = \inf \left\{ \left|   \sum_{i=1}^2 (x_i -  y_i) \theta_i \right| :
\mbx,\mby \in \mcx_{2,2} \text{ and } \mbx \neq \mby  \right\} \\
& =  \inf \left\{ \left|  \theta_1 - \theta_2 \right|, \left| 2 \theta_1 - 2\theta_2 \right|  \right\} 
 =  \left|  \theta_1 - \theta_2 \right| \\
\varepsilon_{2,3}(\Theta_3) 
& = \inf \left\{ \left|   \sum_{i=1}^3 (x_i -  y_i) \theta_i \right| :
\mbx,\mby \in \mcx_{2,3} \text{ and } \mbx \neq \mby  \right\} \\
& =  \inf \left\{ 
\left|  \theta_1 - \theta_2 \right| ,  
\left|  \theta_1 - \theta_3 \right|, 
\left|  \theta_2 - \theta_3 \right|,
\left|  \theta_1 +  \theta_2  - 2\theta_3 \right|, \right. \\
& \hspace{1.5cm}  \left.  \left|  \theta_1 -2 \theta_2 + \theta_3 \right| , 
\left|  -2\theta_1 +  \theta_2 + \theta_3 \right| 
  \right\} \\
  \varepsilon_{2,4}(\Theta_3) 
  & = \inf \left\{ \left|   \sum_{i=1}^4 (x_i -  y_i) \theta_i \right| :
\mbx,\mby \in \mcx_{2,4} \text{ and } \mbx \neq \mby  \right\} \\
& = \inf  \{      \left| \theta_1 - \theta_2\right|, \left| \theta_1 - \theta_3\right|, \left| \theta_1 - \theta_4\right|, \left| \theta_2 - \theta_3\right|, \left| \theta_2 - \theta_4\right|, \left| \theta_3 - \theta_4\right|, \\
&  \qquad \left| 2\theta_1 - \theta_2 - \theta_3\right|, \left| 2\theta_1 - \theta_2 - \theta_4\right|, \left| 2\theta_1 - \theta_3 - \theta_4\right|, \\
& \qquad \left| 2\theta_2 - \theta_3 - \theta_1\right|, \left|  2\theta_2 - \theta_4 - \theta_1\right|, \left| 2\theta_2 - \theta_3 - \theta_4\right|, \\
& \qquad \left|  2\theta_3 - \theta_2 - \theta_1\right|, \left|  2\theta_3 - \theta_2 - \theta_4\right|, \left| 2\theta_3 - \theta_4 - \theta_1\right|, \\
& \qquad  \left|  2\theta_4 - \theta_2 - \theta_1\right|, \left| 2\theta_4 - \theta_2 - \theta_3\right|, \left|  2\theta_4 - \theta_3 - \theta_1\right|, \\
& \qquad  \left| \theta_1 + \theta_2 - \theta_3 - \theta_4 \right|,
 \left| \theta_1 -\theta_2+ \theta_3  - \theta_4 \right|,
  \left| \theta_1 - \theta_2 - \theta_3 + \theta_4\right|     \}. 
\end{align*}

We shall compute  Sidon sets of positive integers.  
Let $h=2$ and $m=1$.  
By Theorem~\ref{Sidon:theorem:construction-integers}, 
in the construction of these sets, we may choose any  integer  
\[
q > \frac{4}{\varepsilon_{2,n}(\Theta_n)}.  
\]
Different choices of $q$ generate different Sidon sets.

Let 
\begin{align*}
\theta_1 = \sqrt{2} & = 1.4142\ldots \\
\theta_2 = \sqrt{3} & = 1.7320\ldots \\
\theta_3 = \sqrt{5} & = 2.2360\ldots \\
\theta_4 = \sqrt{7} & = 2.6457\ldots 
\end{align*} 
and consider the \Q-independent sets 
\[
\Theta_2 = \{\sqrt{2}, \sqrt{3}  \}
\]
\[
\Theta_3 = \{\sqrt{2}, \sqrt{3}, \sqrt{5} \}
\]
\[
\Theta_4 = \{\sqrt{2}, \sqrt{3}, \sqrt{5}, \sqrt{7} \}. 
\]

For $n=2$, we have 
\[
\varepsilon_{2,2}(\Theta_2) = \sqrt{3} - \sqrt{2}  = 0.3178\ldots 
\]
and may choose any integer $q$ such that 
\[
q > \frac{4}{\sqrt{3}-\sqrt{2}} = 12.5850\ldots.
\]
With $q=13$, we have 
\[
\frac{18}{13} < \sqrt{2} < \frac{19}{13} \qqand 
\frac{22}{13} < \sqrt{3} < \frac{23}{13}  
\]
and so 
$A_{2,2}(13,1) = \{a_1,a_2\}$, where 
\begin{align*}
 a_1   & = 18 \text{ or } 19 \\
 a_2  & = 22 \text{ or } 23. 
\end{align*} 
For each of these 4 sets, we have $|2A_{2,2}(13,1)| = 3$.
(Of course, every 2-element set is a Sidon set, so the case $n=2$ is just an exercise in 
working through the \Q-independent construction.)

For $n=3$, we have 
\[
\varepsilon_{2,3}(\Theta_2) =  \sqrt{2} - 2\sqrt{3} + \sqrt{5}  = 0.1861\ldots 
\]
and we choose 
\[
q > \frac{4}{ \sqrt{2} - 2\sqrt{3} + \sqrt{5} } = 21.4845\ldots.
\]
With $q=22$, we have 
\[
\frac{31}{22} < \sqrt{2} < \frac{32}{22}, \hspace{1cm} 
\frac{38}{22} < \sqrt{3} < \frac{39}{22}, \hspace{1cm}   
\frac{51}{22} < \sqrt{5} < \frac{52}{22}  
\]
and so 
$A_{2,3}(22,1) = \{a_1,a_2,a_3\}$, where 
\begin{align*}
 a_1   & = 31 \text{ or } 32 \\
 a_2  & = 38 \text{ or } 39 \\
a_3 & = 51 \text{ or } 52.  
\end{align*} 
For each of these 8 sets we have $|2A_{2,3}(22,1)| = 6$.

For $n=4$, we have 
\[
\varepsilon_{2,4}(\Theta_2) =  \sqrt{2} + \sqrt{7} -  \sqrt{3} - \sqrt{5}  = 0.0918\ldots 
\]
and we choose 
\[
q > \frac{4}{  \sqrt{2} + \sqrt{7} -  \sqrt{3} - \sqrt{5}} = 43.5511\ldots.
\]

With $q=44$, we have 
\[
\frac{62}{44} < \sqrt{2} < \frac{63}{44}, \hspace{1cm} 
\frac{76}{44} < \sqrt{3} < \frac{77}{44},   
\]
\[
\frac{98}{44} < \sqrt{5} < \frac{99}{44}, \hspace{1cm} 
\frac{116}{44} < \sqrt{7} < \frac{117}{44} 
\]
and so 
$A_{2,4}(44,1) = \{a_1,a_2,a_3, a_4\}$, where 
\begin{align*}
 a_1   & =  62\text{ or }  63 \\
 a_2  & =  76 \text{ or } 77\\ 
  a_3  & = 98  \text{ or } 99 \\ 
a_4 & = 116 \text{ or }  117.  
\end{align*} 
For each of these 16 sets, we have $|2A_{2,4}(44,1)| = 10$.

With $q=100$, we have 
\[
\frac{141}{100} < \sqrt{2} < \frac{142}{100}, \hspace{1cm} 
\frac{173}{100} < \sqrt{3} < \frac{174}{100},   
\]
\[
\frac{223}{100} < \sqrt{5} < \frac{224}{100}, \hspace{1cm} 
\frac{264}{100} < \sqrt{3} < \frac{265}{100} 
\]
and so 
$A_{2,4}(100,1) = \{a_1,a_2,a_3, a_4\}$, where 
\begin{align*}
 a_1   & =  141 \text{ or }  142 \\
 a_2  & =  173 \text{ or } 174 \\ 
  a_3  & = 223  \text{ or } 224 \\ 
a_4 & = 264 \text{ or }  265.  
\end{align*} 
For each of these 16 sets, we have $|2A_{2,4}(100,1)| = 10$. 

This example suggests the following result.  

\bt
Let $\Theta_n = \{\theta_1,\ldots, \theta_n\}$ be a \Q-independent set 
of positive real numbers.  For every integer $g \geq 2$, let 
 $\theta_i$ have the $g$-adic expansion 
\[
\theta_i =  \sum_{j=0}^{j_0} c'_{i,j}g^j + \sum_{j=1}^{\infty} \frac{c_{i,j}}{g^j} 
\]
with  digits $c_{i,j} , c'_{i,j} \in [0,g-1]$. 
For every positive integer  $\ell$, let 
\[
a_{i,\ell}  
 = \sum_{j=0}^{j_0} c'_{i,j}g^{\ell +j} + \sum_{j=1}^{\ell}  c_{i,j} g^{\ell-j}. 
\]
There exists $\ell_0 = \ell_0(\Theta_n)$ such that, for all $\ell \geq \ell_0$, the set 
\[
B_{2,n}(\ell) = \{a_{1,\ell},a_{2,\ell},\ldots, a_{n,\ell}\}
\]
is a Sidon set of positive integers. 
\et

Note that $a_{i,\ell}$ is the integer part of $g^{\ell} \theta_i$.

\begin{proof}
For all $i \in [1,n]$, we have 
\[
0 < \theta_i - \frac{a_{i,\ell}}{g^{\ell}} <  \frac{1}{g^{\ell}}.
\]
Apply  Theorem~\ref{Sidon:theorem:construction-integers} with 
\[
q = g^{\ell} > \frac{4}{\varepsilon_{2,n}(\Theta_n)}.
\] 
This completes the proof. 
\end{proof}

\bprob
Let $\Theta_n = \{\theta_1,\ldots, \theta_n\}$ be a set of $n$ positive real numbers 
and let $g,\ell \in \N$ with $g \geq 2$.  
For all $i \in [1,n]$, let $a_{i,\ell}$ be the integer part of $g^{\ell} \theta_i$ 
and let $B_{2,n}(\ell) =  \{a_{1,\ell},a_{2,\ell},\ldots, a_{n,\ell}\}$.  
Suppose that  $B_{2,n}(\ell)$ is a $B_h$-set for some $g \geq 2$ and 
for all $\ell \geq \ell(g)$. 
Does this imply that the set $\Theta_n$ is \Q-independent? 
\eprob

\bprob
Suppose that  $B_{2,n}(\ell)$ is a $B_h$-set for all $g \geq 2$ and 
for all  $\ell \geq \ell(g)$. 
Does this imply that the set $\Theta_n$ is \Q-independent? 
\eprob

\bprob
Let $\Theta$ be an infinite \Q-independent set of real numbers. 
Can one construct from $\Theta$ an infinite Sidon set?
\eprob 

\def\cprime{$'$} \def\cprime{$'$} \def\cprime{$'$}
\providecommand{\bysame}{\leavevmode\hbox to3em{\hrulefill}\thinspace}
\providecommand{\MR}{\relax\ifhmode\unskip\space\fi MR }
\providecommand{\MRhref}[2]{%
  \href{http://www.ams.org/mathscinet-getitem?mr=#1}{#2}
}
\providecommand{\href}[2]{#2}

\end{document}